\documentclass{article}
\usepackage[utf8]{inputenc}
\usepackage{amssymb}
\usepackage{amsmath}
\newtheorem{theorem}{Theorem}

\title{A Combinatorial Proof of a generalization of a Theorem of Frobenius} 
\author{Supravat Sarkar}
\date{}
\begin{document}
\maketitle
\begin{abstract}
    In this article, we shall generalize a theorem due to Frobenius in group theory, which asserts that if $p$ is a prime and $p^{r}$ divides the order of a finite group, then the number of subgroups of order $p^{r}$ is $\equiv$ 1(mod $p$). Interestingly, our proof is purely combinatorial and does not use much group theory.
\end{abstract}
\begin{center}
    \textbf{Keywords}: Group, subgroup, prime, count.
\end{center}
\begin{center}
    \textbf{MSC number}: 20D15
\end{center}

\section*{1.Introduction}

Although Sylow's theorems are taught in almost all undergraduate
courses in abstract algebra, a generalization due to Frobenius does not seem to be as
well known as it ought to be. Frobenius' generalization states that
if $p$ is a prime and $p^{r}$ divides the order $N$ of a finite group
$G$, the number of subgroups of $G$ of order $p^{r}$ is $\equiv 1$
(mod $p$). The special case when $p^{r}$ is the largest power of $p$
dividing $N$ is part of Sylow's third theorem. Many of the standard
texts do not mention this theorem. One source is Ian Macdonald's
`Theory of Groups' \cite{M}. In fact, a further generalization due
to Snapper \cite{S} asserts that for any subgroup $K$ of order
$p^r$ and for any $s \geq r$ where $p^s$ divides the order of $G$,
the number of subgroups of order $p^s$ containing $K$ is also
$\equiv 1$ (mod $p$).  In this article, we give a new proof of a
further extension of Snapper's result that is purely combinatorial
and does not use much group theory. Thus, we have a new
combinatorial proof of Frobenius's theorem as well. \vskip 5mm

\section*{2.Main results}
\vskip 5mm

\noindent We initially started by giving a combinatorial proof of
Frobenius's result and, interestingly, our method of proof yields as
a corollary an extension of Snapper's Theorem. Our proof builds on
the famous combinatorial proof of Cauchy's theorem which asserts
that if a prime divides the order of a group, there is an element of
that prime order. \vskip 5mm

\begin{theorem}
Let $G$ be a finite group of order $N$, and let $p$ be a prime. Let
$b_{0}<b_{1}<\cdots<b_{r}$ be nonnegative integers such that
$p^{b_{r}}$ divides $N$ and $P_{b_{0}}$ be a subgroup of $G$ of
order $p^{b_{0}}$. Then the number of ordered tuples
$(P_{b_{1}},P_{b_{2}},\cdots,P_{b_{r}})$ such that each $P_{b_{i}}$ is
subgroup of $G$ of order $p^{b_{i}}$ and
$$P_{b_{0}}\subset
P_{b_{1}}\subset \cdots \subset P_{b_{r}}$$ is $\equiv 1$ (mod $p$).
\end{theorem}
\vskip 5mm

\noindent  The case $r=1$ is a Theorem due to Snapper \cite{S}
which is itself an extension of Frobenius's Theorem that corresponds
to the case $r=1, b_0=0$ in our Theorem.

\noindent Let us recall here the simple results in finite group
theory that we will need.
\begin{enumerate}
    \item If $H$ is a subgroup of a finite group $G$ of order $N$, and the index $[G:H]$ is the smallest prime divisor of $N$, then $H$ is normal in $G$.
    \item (Sylow's first theorem) If $G$ is a finite group of order $N$, $p$ a prime, $i\geq 0$ is an integer, $p^{i+1}|N$ and $P$ is a subgroup of $G$ of order $p^{i}$, then there is a subgroup $Q$ of $G$ containing $P$ of order $p^{i+1}$.
\end{enumerate}

We shall also use the following notations throughout.
\begin{enumerate}
    \item For a finite set $S$, $|S|$ denotes the number of elements (cardinality) of $S$.
    \item If $G$, $H$ are finite groups, $H\leq G$ means $H$ is a subgroup of $G$.
    \item If $G$, $H$ are finite groups, $H\leq G$, $[G:H]$ denotes the index of $H$ in $G$.
    \item If $G$ is a finite group, the order of $G$ is the number of elements of $G$.
    \item If $H$ is a subgroup of a group $G$, $N_{G}(H)$ denotes the normalizer of $H$ in $G$.
    \item For positive integers $a,b$, we write $a|b$ to mean $a$ divides $b$.
\end{enumerate}
\vskip 3mm

\noindent \textbf{Proof of Theorem.}\\
For ease of understanding, we divide the proof into three steps.
\vskip 3mm

\noindent \textbf{Step 1}: We tackle the case $r=1, b_0=0, b_{1}=1$
first, which just says that if $p$ divides the order of $G$, then
the number of subgroups of $G$ of order $p$ is $\equiv$ 1 (mod $p$).

\noindent Let $T$ = $\{(a_{1},a_{2},...,a_{p})\mid a_{i} \in G$
$\forall i$,
$a_{1}a_{2}...a_{p}=1\}$.\\
Observe that $|T|= N^{p-1}\equiv 0$ (mod $p$), as any choice of
$a_{1},...,a_{p-1}$ uniquely determines $a_{p}.$  Also, if not all
$a_{i}$'s are equal, then $(a_{1},a_{2},...,a_{p})\in T$ implies
$(a_{i},a_{i+1},...,a_{i+p-1})$ for $i=1,2,...,p$ (indices are
modulo $p$) are $p$ distinct elements of $T$. The reason is as
follows:\\
If $(a_{i},a_{i+1},...,a_{i+p-1})= (a_{j},a_{j+1},...,a_{j+p-1})$
for some $i\neq j$, then $a_{k}=a_{k+j-i}\forall k$. By induction,
$a_{k}=a_{k+\alpha(j-i)}$ for any integer $\alpha$. But $i\neq j$
implies $gcd(j-i,p)=1$, as $0<|i-j|<p$ and $p$ is a prime. So, $j-i$
is invertible modulo $p$. So any $1\leq l \leq p$ satisfies $l \equiv
1+\alpha(j-i)$ (mod $p$) for some integer $\alpha$. So,
$a_{l}=a_{1+\alpha(j-i)}=a_{1}$ for any $1\leq l \leq p.$ So,
$a_{l}$'s are all equal, which leads to a contradiction.\\

\noindent So, if $d$ is the number of elements of $G$ of order $p$, then
$0 \equiv |T| \equiv (1+d)$ (mod $p$). So, $d \equiv -1$ (mod $p$) (as there are exactly 1+$d$ elements of $T$ with all $a_{i}$'s equal.) In
each subgroup of order $p$, there are $p-1$ elements of order $p$,
different subgroups of order $p$ intersect at the identity. So,

$-1\equiv d$ =$(p-1)$(number of subgroups of order $p$) $\equiv -$(number of subgroups of order $p$) (mod $p$). 

So, number of subgroups of
order $p$ is $\equiv$ 1 (mod $p$), which finishes the proof for the case $r=1, b_{0}=0, b_{1}=1.$ \vskip 3mm

\noindent \textbf{Step 2:} Now come to a general case. First, we
fix a notation. Let $H$ be any group of order $M$, $p^{n}\vert M,
p^{n+1} \nmid M$, $0\leq r \leq n.$ Let $P_{r}$ be a subgroup of
order $p^{r}$ in $H$. Define
$$S(P_{r},H)=\{(P_{r+1},P_{r+2},\cdots,P_{n}) | P_{i} \leq H, |P_{i}|=p^{i}~\forall~ i,
P_{r}\leq P_{r+1}\leq \cdots \leq P_{n}\leq H\}.$$ So, $S(P_{n},H)$ is a
singleton set, by convention.\\
For $r\leq i<n$ and a subgroup $P_{i}$ of $H$ of order $p^{i}$,
there is a subgroup $P_{i+1}^{\prime}$ of $H$ of order $p^{i+1}$
containing $P_{i}$, by Sylow's theorems.
$[P_{i+1}^{\prime}:P_{i}]=p,$ which is the smallest prime divisor of
$|P_{i+1}^{\prime}|$, so $P_{i}$ is normal in $P_{i+1}^{\prime}$. Hence $P_{i+1}^{\prime}\leq N_{G}(P_{i})$. So, $\frac{P_{i+1}^{\prime}}{P_{i}}$ is a subgroup of order $p$ in $\frac{N_{G}(P_{i})}{P_{i}}$. 
So, $p \vert [N_{G}(P_{i}):P_{i}].$ By the same reasoning, any
subgroup $P_{i+1}$ of $H$ of order $p^{i+1}$ containing $P_{i}$ must
be a subgroup of $N_{G}(P_{i})$, and so $\frac{P_{i+1}}{P_{i}}$ is a subgroup of order $p$ in $\frac{N_{G}(P_{i})}{P_{i}}$. Conversely, any subgroup of order $p$ in $\frac{N_{G}(P_{i})}{P_{i}}$ gives rise via pullback to a subgroup $P_{i+1}$ of $N_{G}(P_{i})$( hence of $H$) of order $p^{i+1}$ containing $P_{i}$. So, there is a one-to-one
correspondence between such $P_{i+1}$ (subgroups of $G$ of order $p^{i+1}$ containing $P_{i}$) and the subgroups of order $p$
of the quotient group $\frac{N_{G}(P_{i})}{P_{i}}$.\\

\noindent So, the number of such $P_{i+1}$ is the number of
subgroups of order $p$ in $\frac{N_{G}(P_{i})}{P_{i}}$, which is
$\equiv$ 1 (mod $p$), in view of Step 1. So, in mod $p$, we can choose $P_{r+1}$ in 1
way, after each such choice we can choose $P_{r+2}$ in 1 way, and so
on. So, $|S(P_{r},H)| \equiv 1$ (mod $p$).\vskip 3mm

\noindent \textbf{Step 3:} Now come to the setup of our theorem. We
have $|S(P_{b_{0}},G)|\equiv 1$(mod $p$), by Step 2. Let us count
$|S(P_{b_{0}},G)|$ in another way. Let $x$ be the number of ordered
tuples as in the statement of our theorem. After choosing any of
such $x$ ordered tuples, we can choose
$(P_{b_{i}+1},...,P_{b_{i+1}-1})$ in $|S(P_{b_{i}},P_{b_{i+1}})|
\equiv$ 1(mod $p$) ways, for each $0\leq i\leq r-1$, and we can choose
$(P_{b_{r}+1},...,P_{n})$ in $|S(P_{b_{r}},G)| \equiv$ 1(mod $p$)
ways.\\
Now, $p^{n}$ is the largest power of $p$ dividing $N$, each $P_{i}$ is a subgroup of $G$ of order $p^{i}$ and
$P_{i} \leq P_{i+1}$ for all $b_{0}\leq i<n$. So, we obtain
$|S(P_{b_{0}},G)| \equiv$ $x$ (mod $p$). Hence finally we get $x
\equiv 1$ (mod $p$) which completes the proof. \vskip 3mm

\noindent \textbf{Remarks} \\The case $r=1$ is Snapper's result and
the further special case $r=1, b_{0}=0$ corresponds to Frobenius'
theorem.

\end{document}